\documentclass[11pt,a4paper]{amsart}
\usepackage[
  a4paper,
  left=3.3cm,
  right=3.3cm,
  top=3.3cm,
  bottom=3.5cm
]{geometry}

\usepackage{amsmath,amssymb,amsthm,mathtools}
\usepackage{graphicx}
\usepackage[hidelinks]{hyperref}

\newtheorem{theorem}{Theorem}[section]
\newtheorem{proposition}[theorem]{Proposition}
\newtheorem{lemma}[theorem]{Lemma}

\theoremstyle{remark}
\newtheorem{remark}[theorem]{Remark}

\newcommand{\Z}{\mathbb Z}

\title[Bak--Sneppen Model on
\(\Z^d\)]{Avalanche Sharpness and the Stationary State of the Bak--Sneppen Model on
\(\Z^d\)}


\author{Olivier Hénard}
\address{Institut de mathématique d'Orsay, Université Paris-Saclay,
91405 Orsay, France}
\email{olivier.henard@universite-paris-saclay.fr}

\begin{document}

\begin{abstract}
We establish the stationary picture for the nearest-neighbour
Bak--Sneppen model on \(\mathbb Z^d\), for every \(d\ge1\).
There exists a critical level \(t_c\) such that, as \(N\to\infty\),
every fixed family of stationary fitnesses on the torus
\((\mathbb Z/N\mathbb Z)^d\) converges, in exponential coordinates,
to independent random variables with law
\(t_c+\operatorname{Exp}(1)\).

The proof relies on three estimates for avalanches. For \(d\ge2\), these
replace the interval structure available in dimension one. We show that
the threshold at which the mean avalanche range diverges coincides with
the threshold for the appearance of an infinite avalanche, using the
hierarchical structure of avalanches and an estimate on the creation of
new sites under translation. We then prove that an avalanche reaching a
macroscopic distance on the torus can, after any fixed positive increase
of the level, be amplified into one covering the whole torus. Finally,
below criticality, we obtain  exponential integrability of
the avalanche range and give a finite--infinite volume comparison.

These three ingredients allow us to implement the Meester and
Znamenski program that connects these avalanche estimates to the
stationary law of the finite-volume process.
\end{abstract}

\maketitle

\section{Introduction and main results}
\label{sec:introduction}

The Bak--Sneppen model is a simple model of self-organized criticality
\cite{BakSneppen1993}. Each vertex of a finite graph carries a fitness,
and at every step the vertex of smallest fitness and its neighbours
receive fresh independent values. Despite the simplicity of this rule,
the stationary fitness distribution is expected to develop, as the
system grows, a sharp threshold: fitnesses typically lie above this
level, while the dynamics produces bursts of updates involving smaller
values. These bursts are called avalanches. Fixing a level \(t\) and
following the updates until all fitnesses are again above \(t\)
provides a natural way to isolate and study them
\cite{Maslov1996,PaczuskiMaslovBak1996}.

Meester and Znamenski \cite{MeesterZnamenski2003,MeesterZnamenski2004}
proposed a rigorous program connecting these avalanches to the stationary
large-volume limit; see also
\cite{GillettMeesterNuyens2006,GillettMeesterVanderWal2006,PopovVolkov2026,Volkov2022}.
The key is to show that three critical levels coincide: the level at
which the mean avalanche duration becomes infinite, the level at which
the mean number of sites visited becomes infinite, and the level at
which an avalanche can last forever. Together with suitable
finite-volume estimates, this identifies the stationary limit as a
product measure supported above a common fitness threshold. We recall
this program in Section~\ref{sec:stationary}.

\begin{figure}[t]
\centering
\includegraphics[width=0.66\textwidth]{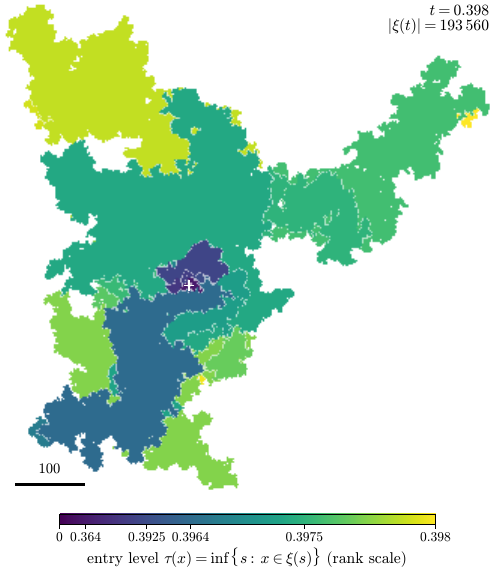}
\caption{A near-critical avalanche on $\Z^2$ at level
$t=0.398$ ($t_c\approx0.3988$).
Sites are coloured by their entry level $\tau(x)$ in the
monotone coupling, revealing the nested avalanche ranges.
Of the 133 grafted avalanches, only 43 reach new sites,
accounting for the 44 colours.
The $+$ marks the origin.
This realization is the largest among approximately
$9000$ independent avalanches at this level. (Most avalanches are still small at criticality).}
\label{fig:avalanche}
\end{figure}

For the one-dimensional one-sided model, the range size admits a closed
scalar evolution \cite{HenardSharpness}. In the isotropic one-dimensional
model, the range remains an interval, but the scalar description based
on the graft length alone no longer closes, since the locations of the two
endpoints matter. For \(d\ge2\), avalanche geometry is even more complex,
largely unknown, and no comparable scalar description is available.
Figure~\ref{fig:avalanche} illustrates this nested geometry in
dimension two. Our aim is to replace these one-dimensional reductions
by estimates that remain effective in arbitrary dimension. A key
quantity is the overlap of independent avalanche ranges.
% ; the main point is that an upper bound on its
% expectation already suffices for sharpness.

% Simulation: ./bs aval 1024 0.398 100000000 8299 avb8299 ; python3 make_figure_avalanche.py
% t_c ~ 0.3988 = -log(1-f_c), f_c ~ 0.3289 (2D estimate, Paczuski--Maslov--Bak 1996 -- to be checked).

We consider the model on the discrete torus
\[
\Lambda_N=(\Z/N\Z)^d,\qquad d\ge1,
\]
and use exponential coordinates for the fitnesses. Thus every
refreshed fitness has law \(\operatorname{Exp}(1)\), and an update at
\(x\) refreshes the sites in
\[
x+\mathcal N_d,
\qquad
\mathcal N_d:=\{0,\pm e_1,\ldots,\pm e_d\}.
\]
Let \(\pi_N\) be the stationary law of this finite-volume Markov chain, and, for
\(x\in\Z^d\), let \(F_{N,x}\) denote the stationary fitness at the image of
\(x\) in \(\Lambda_N\).

Our main result identifies the large-volume stationary state.

\begin{theorem}[Stationary limit]
\label{thm:stationary-main}
There exists \(t_c\in(0,\infty)\) such that, for every fixed collection of
distinct sites \(x_1,\ldots,x_k\in\Z^d\),
\[
(F_{N,x_1},\ldots,F_{N,x_k})
\xrightarrow[N\to\infty]{\mathrm d}
(t_c+E_1,\ldots,t_c+E_k),
\]
where \(E_1,\ldots,E_k\) are independent
\(\operatorname{Exp}(1)\) random variables.
\end{theorem}

For \(y\ge0\), let \(\nu_y\) denote the law of \(y+E\), where
\(E\sim\operatorname{Exp}(1)\). Under the increasing change of variables
\[
u=1-e^{-f},
\]
Theorem~\ref{thm:stationary-main} is equivalent to the usual uniform
representation: the limiting fitnesses are independent and uniform on
\[
[1-e^{-t_c},1].
\]

The critical level \(t_c\) is characterized through infinite-volume
avalanches. Fix \(t\ge0\), start from a configuration in which
all fitnesses are larger than \(t\), and force an update at the origin.
Thereafter, as long as some current fitness is at most \(t\), update the
site carrying the smallest fitness and its neighbours. The procedure
stops when all current fitnesses are again larger than \(t\). The set of
sites refreshed during this procedure (including the forced initial update) is the range of the
\(t\)-avalanche.

Let \(\xi(t)\) denote this range on \(\Z^d\), let \(\eta(t)\) be its
duration, including the forced update and with value \(\infty\) if the
avalanche does not terminate, and put
\[
r(t):=|\xi(t)|,\qquad R(t):=\mathbb E r(t),\qquad
D(t):=\mathbb E\eta(t),
\qquad \theta(t):=\mathbb P(r(t)=\infty).
\]
We use throughout the natural monotone coupling in the level parameter,
under which \(\xi(s)\subseteq\xi(t)\) for \(s\le t\). Avalanche ranges are
connected, since every newly refreshed site is a nearest neighbour of a
site already reached by the avalanche.

Following Meester and Znamenski, introduce the three critical levels
\[
\begin{aligned}
t_{\rm d}&:=\sup\{t:D(t)<\infty\},
& t_{\rm sus}&:=\sup\{t:R(t)<\infty\},
& t_\infty&:=\inf\{t:\theta(t)>0\}.
\end{aligned}
\]
Every update refreshes at most \(2d+1\) sites, hence
\(r(t)\le(2d+1)\eta(t)\), while positive mass at infinity forces infinite
mean range. Therefore
\[
t_{\rm d}\le t_{\rm sus}\le t_\infty.
\]
Both inequalities could a priori be strict. Our first result identifies all
three levels.

\begin{theorem}[Avalanche sharpness]
\label{thm:sharpness-main}
For the nearest-neighbour Bak--Sneppen model on \(\Z^d\), \(d\ge1\),
\[
t_{\rm d}=t_{\rm sus}=t_\infty=:t_c.
\]
Moreover, for every \(t<t_c\),
\begin{equation}
\label{eq:intro-sus-bound}
t_c  \le t + \frac{10}{R(t)}.
\end{equation}
\end{theorem}

\begin{remark}
The equality
\[
t_{\rm d}=t_{\rm sus}
\]
was already derived by Meester and Znamenski \cite[Theorem~3.3 and Section~6]{MeesterZnamenski2004}, through the
range--duration relation, stated here as Lemma~\ref{lem:range-duration}. The novelty here is the
identification \(t_{\rm sus}=t_\infty\), together with the quantitative
bound \eqref{eq:intro-sus-bound}.
\end{remark}
\begin{remark}
Finiteness of $t_c$ follows
directly from \eqref{eq:intro-sus-bound}, taking \(t=0\), since $R(0)=2d+1$.
Positivity of \(t_c\) is elementary:  in the exploration of a
\(t\)-avalanche, each update produces \(2d+1\) fresh fitnesses, each of
which is at most \(t\) with probability \(1-e^{-t}\). Ignoring overlaps
therefore gives a domination by a Galton--Watson process with offspring
distribution
\[
\operatorname{Bin}(2d+1,1-e^{-t}).
\]
This process is subcritical whenever
$(2d+1)(1-e^{-t})<1$, which implies finite mean avalanche duration for all sufficiently small
\(t>0\), and hence \(t_c>0\).
\end{remark}

For \(d=1\), Theorem~\ref{thm:sharpness-main} already suffices to complete
the Meester--Znamenski program \cite{MeesterZnamenski2003}, recalled in
Section~\ref{sec:stationary}, and derive the stationary-limit result. In
dimensions \(d\ge2\), two additional finite-volume ingredients are
required.\\

The first ingredient concerns avalanches on the finite torus. Let
\(\xi_N(t)\) denote the range of the \(t\)-avalanche on \(\Lambda_N\), and
set
\[
r_N(t):=|\xi_N(t)|,\qquad
R_N(t):=\mathbb E r_N(t),\qquad
P_N(t):=\mathbb P(\xi_N(t)=\Lambda_N).
\]
The following result says that once a supercritical avalanche has reached
a macroscopic distance, any further fixed increase of the level amplifies
it into complete coverage of the torus.

\begin{theorem}[Macroscopic amplification]
\label{thm:amplification-main}
Fix
\[
t_c<s<t
\qquad\text{and}\qquad
0<\rho<\frac12.
\]
Let
\[
E_N(s)
:=
\left\{
\xi_N(s)\cap
\{x:d_{\Lambda_N}(0,x)=\lfloor\rho N\rfloor\}
\neq\varnothing
\right\}.
\]
Then
\[
\mathbb P\bigl(
\xi_N(t)=\Lambda_N
\mid E_N(s)
\bigr)
\longrightarrow1.
\]
Consequently,
\[
\liminf_{N\to\infty}P_N(t)
\ge
\theta(s)>0.
\]
\end{theorem}

The complementary statement below shows that, in the subcritical regime,
an avalanche is typically too small to detect the periodic geometry of
the torus. 

\begin{theorem}[Subcritical exponential moments and locality]
\label{thm:subcritical-main}
For every \(t<t_c\), there exists \(\lambda_t>0\) such that
\[
\mathbb E[e^{\lambda_t r(t)}]<\infty.
\]
In particular, all polynomial moments of \(r(t)\) are finite. Moreover,
\[
\sup_{0\le u\le t}|R_N(u)-R(u)|
\longrightarrow0.
\]
\end{theorem}

% \begin{remark}
% In dimension \(d=1\), the two additional ingredients just described are
% either trivial or not needed to establish the stationary limit result:
% macroscopic amplification is immediate from the interval structure of
% avalanche ranges, while the required subcritical finite-volume
% approximation follows from a finite first moment alone.
% \end{remark}

\subsection{Proof strategy}

The proof for general \(d\ge1\) is organized around the following three
mechanisms.

\medskip
\noindent
\textbf{1. Infinite-volume sharpness: creation under translation.}
The level parameter gives avalanches a hierarchical structure.
Conditionally on a finite current range \(A\), each site \(x\in A\)
produces marks at rate one, each grafting an independent avalanche
translated by \(x\). The relevant geometric quantity is the part of this
translated range lying outside \(A\).

The deterministic input is
\[
\sum_{x\in A}|(x+B)\setminus A|
\ge c|A||B|,
\qquad 0\in B,\quad |B|\ge|A|\ge3,
\]
which follows from a translation-invariant total order on \(\Z^d\).
Applied to
\[
\psi(t,\lambda):=\mathbb E[1-e^{-\lambda r(t)}],
\]
the hierarchical decomposition gives the Burgers-type inequality
\[
\partial_t\psi(t,\lambda)
\ge
\frac1{10}\psi(t,\lambda)\partial_\lambda\psi(t,\lambda).
\]
Since \(\psi(t,0+)=\theta(t)\), while
\(\partial_\lambda\psi(t,0+)=R(t)\) whenever \(r(t)<\infty\) almost surely,
the characteristic flow prevents the boundary slope from becoming
infinite strictly before positive boundary mass appears. This yields
Theorem~\ref{thm:sharpness-main}.

\medskip
\noindent
\textbf{2. Supercritical amplification: translation expansion.}
An avalanche reaching distance of order \(N\) need not resemble a ball or
be well mixed on the torus. We retain only one nearest-neighbour edge
crossing each successive radius.

For \(A\subset\Lambda_N\), let
\[
b_A(z):=|\{x\in A:x+z\notin A\}|
\]
be its translation boundary. It satisfies
\[
b_A(z+w)\le b_A(z)+b_A(w),
\qquad
\sum_{z\in\Lambda_N}b_A(z)=|A|(N^d-|A|).
\]
Fix \(0<\rho<1/2\), and let \(E\) be the event that an independent
avalanche reaches distance \(\rho N\). On \(E\), select uniformly one
crossing edge at each radius \(1\le j\le L:=\lfloor\rho N\rfloor\).
Let \(\kappa\) be their expected endpoint counting measure, with zero
contribution off \(E\). The selected edges need not form a connected path.

At a fixed supercritical level, \(\mathbb P(E)\) is bounded away from zero
uniformly in \(N\), so \(\kappa(\Lambda_N)=2L\mathbb P(E)\) is of order
\(N\). Hyperoctahedral symmetry, subadditivity of \(b_A\), and the crossing
edges then imply
\[
\sum_z\kappa(z)b_A(z)
\ge
cN\,\frac{|A|(N^d-|A|)}{N^d}.
\]
Thus the averaged endpoint measure supplies enough translations to expand
any nontrivial set \(A\).

To turn expansion into covering, we use the Lyapunov function
\[
\Phi(A):=\frac{N^d}{|A|}-1.
\]
Conditionally on macroscopic reach at the initial level, the preceding
estimate implies exponential decay of the mean of \(\Phi\) at rate of
order \(N\). Since every noncovering range satisfies
\(\Phi\ge1/(N^d-1)\), any fixed positive sprinkling then forces complete
coverage with high (conditional) probability.

\medskip
\noindent
\textbf{3. Subcritical locality: exponential moments replace exact geometry.}
Below \(t_c\), the hierarchy gives a transport inequality for the
truncated exponential transform of the avalanche range. Finite
susceptibility can thus be bootstrapped to exponential integrability of
\(r(t)\), making ranges of order \(N\) exponentially unlikely.

Every connected subset of \(\Z^d\) with fewer than \(N\) vertices projects
injectively onto \(\Lambda_N\). A subcritical avalanche is therefore
overwhelmingly unlikely to detect the periodic geometry. This gives the
finite--infinite volume comparison of Theorem~\ref{thm:subcritical-main}.

\medskip

These mechanisms provide the inputs for the locking-threshold
representation; see \eqref{eq:locking-recursion} and
\cite[Section~2]{MeesterZnamenski2003}. Sharpness identifies the candidate
threshold, while amplification and subcritical locality control the
locking thresholds from above and below. The stationary limit then
follows from a dimension-free finite-volume argument.

\medskip
\noindent
\textit{Choice of the underlying graph.}
Several parts of the argument extend beyond \(\Z^d\), for instance
to suitable Cayley graphs and their finite quotients. One could formulate
the proof in such an abstract setting, replacing translations and the
amplification argument by appropriate group-theoretic and spectral
assumptions. We have chosen not to do so: the lattice \(\Z^d\)
already contains the essential difficulties, and its geometry makes the
amplification mechanism particularly transparent.

\subsection{Methodological perspective}

Two analogies help place the proof. The first is with coagulation and
gelation. For solvable Smoluchowski coagulation equations, Laplace
transforms lead to first-order nonlinear equations related to inviscid
Burgers, while gelation corresponds to mass escaping from the finite-cluster
distribution; see \cite{Aldous1999,MenonPego2004}. Spatial overlaps prevent
a closed equation for the avalanche law, but the translation estimate still
yields
\[
\partial_t\psi
\ge c\,\psi\,\partial_\lambda\psi.
\]
This one-sided inequality is enough for the characteristic argument:
divergence of the susceptibility is transported into positive boundary
mass at a higher level, signalling an infinite avalanche.

The second analogy concerns amplification. The measure \(\kappa\) may be
viewed as the jump measure of a translation walk on \(\Lambda_N\), and the
expansion estimate is equivalent, up to constants, to a spectral-gap bound
for this auxiliary walk; see Remark~\ref{rem:poincare} and compare
\cite{DiaconisSaloffCoste1993}. Here the estimate follows directly from
subadditivity of the translation boundary and the selected crossing edges,
without any shape or mixing assumption on the macroscopic range.

\subsection*{Disclosure on the use of AI}
Human was used at the level of an annoyingly slow co-author in this work (see Appendix~\ref{appendix-AI} for a more detailed account of our interactions with AI). In view of the extent of the AI contribution, we do not presently intend to submit this work for publication.

\medskip
\noindent
\textit{Relation with \cite{HenardSharpness}.}
The present paper supersedes the first version of \cite{HenardSharpness},
whose main goal was to establish sharpness in dimension 1. Accordingly, that work
has been refocused on quantitative bounds for the \textit{critical tail} of the
avalanche range for the anisotropic model. This is a genuinely different question, on which the
present method, despite its strength for proving sharpness, has little to say.

\subsection{Organization of the paper}

Section~\ref{sec:hierarchy} gives the avalanche hierarchy and its basic
identities.
Section~\ref{sec:burgers-comparison} proves infinite-volume sharpness.
Section~\ref{sec:amplification} establishes the finite-volume expansion
inequality and the amplification theorem.
Section~\ref{sec:subcritical-moments} proves exponential moments of the
subcritical range and the finite--infinite volume comparison.
The stationary limit is deduced in Section~\ref{sec:stationary}.

\section{Avalanche hierarchy and basic identities}
\label{sec:hierarchy}

We collect here the structural facts about avalanches which are common to
all three parts of the proof.

Recall that \(\xi(t)\) and \(\xi_N(t)\) denote the ranges of the
\(t\)-avalanches rooted at the origin in \(\Z^d\) and \(\Lambda_N\),
respectively. We use the standard monotone
graphical construction in the level parameter, so that
\[
\xi(s)\subseteq\xi(t),
\qquad s\le t.
\]

The key property is regeneration in the level variable. We use the
censored exploration filtration \((\mathcal H_u)_{u\ge0}\): during the
exploration of a \(u\)-avalanche, refreshed values at most \(u\) are
revealed exactly, whereas values larger than \(u\) are revealed only
through the event that they exceed \(u\). The \(\sigma\)-field
\(\mathcal H_u\) records the explored genealogy, the identities of the
refresh variables encountered, and this partially revealed information.
The monotone graphical construction makes this filtration increasing.

\begin{proposition}[Hierarchical avalanche decomposition]
\label{prop:hierarchical-decomposition}
Consider either the infinite lattice or the finite torus.
On the event that the range \(A\) at level \(u\) is finite, conditionally
on \(\mathcal H_u\), the terminal fitnesses at its sites are independent
with law \(\nu_u\). Each site \(x\in A\) therefore carries an independent
rate-one clock in the level parameter. At a mark rooted at \(x\) and
level \(u\), an independent \(u\)-avalanche \(\xi'(u)\) is grafted onto
the current range:
\[
A\longmapsto A\cup(x+\xi'(u)).
\]
The graft is independent of the censored past. In infinite volume this
description is used on finite-range states, up to explosion, or with the
process stopped when the range reaches a prescribed finite size.
\end{proposition}

\begin{proof}
We first verify that finite range implies finite duration almost surely, so that the terminal fitness at every reached site is well defined.
Indeed, while the range has at most \(K\) sites, there are at most \(K\)
active fitnesses at or below \(u\). A block of \(K\) updates in which all
fresh values exceed \(u\) removes these active fitnesses without creating
new ones, and therefore forces termination within the block. Such a
block has conditional probability at least
\[
e^{-(2d+1)uK}>0.
\]
The probability of completing \(m\) blocks without either terminating or
reaching more than \(K\) sites is consequently at most
\((1-e^{-(2d+1)uK})^m\). Letting \(m\to\infty\), and then taking the union
over \(K\), proves the assertion.

Since the avalanche has finite duration, each reached site has a last refresh. The censored exploration determines the identities of these last refresh variables, without revealing their values beyond the fact that they exceed \(u\). As these variables are distinct, independence and the memoryless property imply that, conditionally on \(\mathcal H_u\), the terminal fitnesses are independent with law \(u+\operatorname{Exp}(1)\).

% The censored exploration of this finite avalanche identifies the last
% refresh variable assigned to each reached site. These variables are distinct, and each has been observed only to
% exceed \(u\). Independence and the lack-of-memory property give independent
% terminal values \(u+\operatorname{Exp}(1)\), conditionally on
% \(\mathcal H_u\). 
% Their excesses supply the independent level clocks.
At the next triggering level, fresh refresh variables generate an
independent avalanche at that level, translated to the triggering site.
Iterating this argument gives the stated grafting description.

Equivalently, on the torus over any finite level interval, and in infinite
volume over \([0,T]\) with \(T<t_\infty\), independently for each site
\(x\), take Poisson marks \((u,B)\) with intensity
\(du\,\mathbb P(\xi_N(u)\in dB)\) on the torus and
\(du\,\mathbb P(\xi(u)\in dB)\) in infinite volume. Starting from
\(\Xi_0=\mathcal N_d\), retain a mark \((u,B)\) whenever \(x\) belongs to
the current range, in which case
\[
\Xi_u=\Xi_{u-}\cup(x+B).
\]
This realizes the preceding level-by-level description. On the torus it
is directly non-explosive. In infinite volume, monotonicity gives
\[
r(u)
\le
r(T)<\infty
\qquad\text{a.s.},\qquad u\le T,
\]
and one may equivalently stop when the range reaches \(K\) sites and then
let \(K\to\infty\).
\end{proof}

Thus, conditionally on a current range \(A\), marks occur at total rate
\(|A|\). This observation is the common source of the compensation identity
of Section~\ref{sec:burgers-comparison}, the growth estimates of
Section~\ref{sec:amplification}, and the transport inequality of
Section~\ref{sec:subcritical-moments}. We shall also use the same
construction with each graft decorated by the duration of the corresponding
independent avalanche; this gives the range--duration identity below.

Let \(\eta_N(t)\) denote the duration of the finite-volume
\(t\)-avalanche. With \(\eta(t)\) and \(D(t)\) as above, set
\[
D_N(t):=\mathbb E\eta_N(t).
\]

\begin{lemma}[Range--duration identity]
\label{lem:range-duration}
For every \(t<t_{\rm sus}\), \(D(t)<\infty\), and on the torus
\(D_N(t)<\infty\) for every finite \(t\). Moreover, \(R\) is continuous on
\([0,t_{\rm sus})\), \(R_N\) is continuous on \([0,\infty)\), and \(D\) and
\(D_N\) are \(C^1\) on these respective domains, with
\[
D'(t)=R(t)D(t),
\qquad
D_N'(t)=R_N(t)D_N(t),
\]
and
\[
D(0)=D_N(0)=1.
\]
Consequently,
\begin{equation}
\label{eq:range-duration}
D(t)
=
\exp\left\{\int_0^tR(u)\,du\right\},
\qquad
D_N(t)
=
\exp\left\{\int_0^tR_N(u)\,du\right\}.
\end{equation}
\end{lemma}

\begin{proof}
Fix \(T<t_{\rm sus}\), so that \(R(u)\le R(T)<\infty\) for \(u\le T\).
In the decorated hierarchical construction a mark at level \(u\) adds
exactly the duration \(\eta'(u)\) of an independent \(u\)-avalanche.
For \(L\ge1\), put \(D_L(u):=\mathbb E[\eta(u)\wedge L]\). Since
\[
(\eta+\eta')\wedge L-\eta\wedge L\le\eta'\wedge L,
\]
compensation, first stopped when the range reaches \(K\) sites, gives
after letting \(K\to\infty\)
\[
D_L(t)\le1+\int_0^tR(u)D_L(u)\,du,
\qquad 0\le t\le T.
\]
Here the rate of marks is the current range size, and the graft is
independent of the censored past. The stopping can be removed because
\(r(T)<\infty\) almost surely, the observable is bounded by \(L\), and
its expected compensator is bounded by \(LR(u)\). Gronwall's inequality
yields
\[
D_L(t)\le\exp\left\{\int_0^tR(u)\,du\right\}.
\]
Monotone convergence as \(L\to\infty\) proves \(D(t)<\infty\), uniformly
on \([0,T]\). The same argument on the torus uses \(R_N(u)\le N^d\) and
proves \(D_N(t)\le e^{N^dt}<\infty\) on every bounded level interval.

We next record the continuity of the mean ranges. Couple the avalanches
\(\xi(u)\), \(0\le u\le T\), through the same hierarchical construction.
Since \(r(T)<\infty\) almost surely, only finitely many relevant marks occur
below level \(T\), and no deterministic level carries a mark. Hence, for
any deterministic \(t_n\to t\) in \([0,T]\), one has
\(\xi(t_n)=\xi(t)\), and therefore \(r(t_n)=r(t)\), eventually almost
surely. Since \(r(t_n)\le r(T)\) and
\(\mathbb E r(T)=R(T)<\infty\), dominated convergence yields
\(R(t_n)\to R(t)\). Thus \(R\) is continuous on \([0,t_{\rm sus})\).
The same argument on the torus is immediate from \(r_N(T)\le N^d\), and
shows that \(R_N\) is continuous on every bounded level interval.

We may now apply compensation to the full duration. Tonelli's theorem
gives
\[
D(t)=1+\int_0^tR(u)D(u)\,du,
\qquad
D_N(t)=1+\int_0^tR_N(u)D_N(u)\,du.
\]
In infinite volume the stopping by range size is again removed by
monotone convergence of the accumulated duration and compensator.
The integral identities first give continuity of \(D\) and \(D_N\).
Together with the continuity of \(R\) and \(R_N\), the fundamental theorem
of calculus then yields the stated \(C^1\) identities. Since
\(D,D_N\ge1\), integrating \((\log D)'=R\) and
\((\log D_N)'=R_N\) gives \eqref{eq:range-duration}.
\end{proof}

\section{Infinite-volume sharpness}
\label{sec:burgers-comparison}

We prove Theorem~\ref{thm:sharpness-main}. The hierarchical avalanche
decomposition provides an exact compensation formula for a bounded
transform of the range, while a deterministic translation estimate ensures
that an independent graft creates sufficiently many new sites on average.
Their combination yields a Burgers-type differential inequality.

For \(\lambda>0\), put
\[
\psi(t,\lambda)
:=
\mathbb E[1-e^{-\lambda r(t)}]
\]
with the convention
\[
e^{-\lambda\infty}=0.
\]

\begin{lemma}[Avalanche compensation]
\label{lem:avalanche-compensation}
For every \(\lambda>0\) and \(0\le s<t<t_\infty\),
\begin{equation}
\label{eq:compensation-integrated}
\psi(t,\lambda)-\psi(s,\lambda)
=
\int_s^tH(u,\lambda)\,du,
\end{equation}
where, if \(\xi'(u)\) is an independent copy of
\(\xi(u)\),
\begin{equation}
\label{eq:compensator}
H(u,\lambda)
:=
\mathbb E\left[
e^{-\lambda r(u)}
\sum_{x\in\xi(u)}
\left(
1-
e^{-\lambda|(x+\xi'(u))\setminus\xi(u)|}
\right)
\right].
\end{equation}
In particular, \(t\mapsto\psi(t,\lambda)\) is locally Lipschitz on
\([0,t_\infty)\). The functions \(H\) and \(\partial_\lambda\psi\) are
jointly continuous on \([0,t_\infty)\times(0,\infty)\), and
\[
\partial_t\psi(t,\lambda)=H(t,\lambda),
\qquad
\partial_\lambda\psi(t,\lambda)
=\mathbb E[r(t)e^{-\lambda r(t)}].
\]
The \(t\)-derivative at \(t=0\) is understood from the right.
\end{lemma}

\begin{proof}
Fix \(T<t_\infty\), and stop the hierarchical construction when the range
reaches \(K\) sites. Before this stopping time the total rate of relevant
marks is at most \(K\), so compensation applies without any integrability
assumption.

For
\[
F_\lambda(A):=1-e^{-\lambda|A|},
\]
a graft \(B\) rooted at \(x\in A\) produces the increment
\[
F_\lambda(A\cup(x+B))-F_\lambda(A)
=
e^{-\lambda|A|}
\left(
1-e^{-\lambda|(x+B)\setminus A|}
\right).
\]
Compensation for the stopped construction therefore gives the stopped
version of \eqref{eq:compensation-integrated}.

The compensator is uniformly bounded in the stopping level, since
\[
|A|e^{-\lambda|A|}
\le\frac1{e\lambda}.
\]
Moreover,
\[
\mathbb P(r(T)\ge K)\longrightarrow0.
\]
Dominated convergence as \(K\to\infty\) proves
\eqref{eq:compensation-integrated} and \eqref{eq:compensator}. The same
bound gives
\[
|\psi(t,\lambda)-\psi(s,\lambda)|
\le
\frac{|t-s|}{e\lambda},
\]
hence local Lipschitz continuity.

For joint continuity, couple two independent range processes on
\([0,T]\), with \(T<t_\infty\). The pathwise local-constancy argument established
in the proof of Lemma~\ref{lem:range-duration} applies verbatim: since
\(T<t_\infty\), the terminal ranges are finite almost surely, and for any
deterministic \(t_n\to t\) in \([0,T]\), both ranges are eventually
unchanged almost surely. No moment assumption is needed here. If also
\(\lambda_n\to\lambda>0\), the integrand defining \(H\) converges almost
surely and is bounded by \(1/(e\lambda_*)\), where \(\lambda_*>0\) is a
common lower bound for the exponents. Dominated convergence gives joint
continuity of \(H\). The same bound justifies differentiation of \(\psi\)
under the expectation and joint continuity of \(\partial_\lambda\psi\);
the integrand of \(\psi\) itself is bounded by one. Finally,
\eqref{eq:compensation-integrated} and continuity of \(H\) yield
\(\partial_t\psi=H\) at every level in the stated domain.
\end{proof}

The geometric input is the following elementary inequality.

\begin{lemma}[Average creation under translation]
\label{lem:translation-creation}
Let \(A,B\subset\Z^d\) be finite, with \(0\in B\), and write
\[
a:=|A|,
\qquad
b:=|B|.
\]
If \(b\ge a\ge3\), then
\begin{equation}
\label{eq:translation-creation}
\sum_{x\in A}|(x+B)\setminus A|
\ge\frac15ab.
\end{equation}
\end{lemma}

\begin{proof}
Write \(C(A,B):=\sum_{x\in A}|(x+B)\setminus A|\), and first assume
\(b=a\). Equip \(\Z^d\) with the translation-invariant lexicographic
order, and write
\[
A=\{a_1<\cdots<a_a\},
\qquad
p:=|\{z\in B:z>0\}|,
\qquad
q:=|\{z\in B:z<0\}|.
\]
Then \(p+q=a-1\). From \(a_i\), at most \(a-i\) positive displacements
and \(i-1\) negative displacements can land in \(A\). Hence
\[
C(A,B)\ge
\sum_{i=1}^a\bigl((p-a+i)_++(q-i+1)_+\bigr)
=\frac{p(p+1)+q(q+1)}2
\ge\left\lfloor\frac{a^2}4\right\rfloor.
\]
Here \(v_+:=\max\{v,0\}\).

For general \(b\ge a\), choose uniformly a subset \(B_0\subset B\) of
size \(a\) containing zero. The zero displacement contributes nothing,
while every other displacement is retained with probability
\((a-1)/(b-1)\). Applying the equal-cardinality case gives
\[
\frac{a-1}{b-1}C(A,B)
=\mathbb E[C(A,B_0)]
\ge\left\lfloor\frac{a^2}4\right\rfloor.
\]
Since \((b-1)/(a-1)\ge b/a\) and \(a\ge3\),
\[
C(A,B)\ge
\frac ba\left\lfloor\frac{a^2}4\right\rfloor
\ge\frac15ab,
\]
as claimed.
\end{proof}

\begin{remark}
\label{rem:ordered-group}
The proof of Lemma~\ref{lem:translation-creation} uses only the existence
of a translation-invariant total order and the cardinalities of \(A\) and
\(B\). In particular, neither isoperimetry nor any metric property of
\(\Z^d\) enters the infinite-volume sharpness argument.
\end{remark}

We now combine the two lemmas.

\begin{proposition}[Burgers comparison]
\label{prop:burgers-comparison}
For every \(0\le s<t<t_\infty\) and \(\lambda>0\),
\begin{equation}
\label{eq:burgers-integrated}
\psi(t,\lambda)-\psi(s,\lambda)
\ge
\frac1{10}
\int_s^t
\psi(u,\lambda)\,
\partial_\lambda\psi(u,\lambda)\,du.
\end{equation}
Equivalently,
\[
\partial_t\psi(t,\lambda)
\ge
\frac1{10}
\psi(t,\lambda)\partial_\lambda\psi(t,\lambda)
\]
for every \(0\le t<t_\infty\).
\end{proposition}

\begin{proof}
Fix \(u<t_\infty\), and let
\[
A:=\xi(u),
\qquad
B:=\xi'(u),
\qquad
X:=|A|,
\qquad
Y:=|B|.
\]
Then \(X\) and \(Y\) are independent and identically distributed, finite
almost surely, and
\[
X,Y\ge2d+1\ge3.
\]
For \(x\in A\), put
\[
d_x:=|(x+B)\setminus A|.
\]
On \(\{Y\ge X\}\), Lemma~\ref{lem:translation-creation} gives
\[
\sum_{x\in A}d_x
\ge
\frac15XY.
\]

Since \(r\mapsto1-e^{-\lambda r}\) is concave, vanishes at zero, and
\(0\le d_x\le Y\),
\[
1-e^{-\lambda d_x}
\ge
\frac{d_x}{Y}(1-e^{-\lambda Y}).
\]
Consequently, on \(\{Y\ge X\}\),
\[
\sum_{x\in A}
(1-e^{-\lambda d_x})
\ge
\frac X5(1-e^{-\lambda Y}).
\]
Lemma~\ref{lem:avalanche-compensation} therefore gives
\begin{equation}
\label{eq:pre-burgers}
H(u,\lambda)
\ge
\frac15
\mathbb E\left[
Xe^{-\lambda X}
(1-e^{-\lambda Y})
\mathbf1_{\{Y\ge X\}}
\right].
\end{equation}

Set
\[
f(r):=re^{-\lambda r},
\qquad
g(r):=1-e^{-\lambda r}.
\]
The ratio
\[
\frac{f(r)}{g(r)}
=
\frac{r}{e^{\lambda r}-1}
\]
is decreasing on \((0,\infty)\). Hence, on \(\{X<Y\}\),
\[
f(Y)g(X)\le f(X)g(Y).
\]
By exchangeability,
\[
\begin{aligned}
\mathbb E[
f(X)g(Y)\mathbf1_{\{Y<X\}}
]
&=
\mathbb E[
f(Y)g(X)\mathbf1_{\{X<Y\}}
] \\
&\le
\mathbb E[
f(X)g(Y)\mathbf1_{\{X<Y\}}
].
\end{aligned}
\]
Therefore
\[
\mathbb E[
f(X)g(Y)\mathbf1_{\{Y\ge X\}}
]
\ge
\frac12
\mathbb E[f(X)g(Y)].
\]
Independence gives
\[
\mathbb E[f(X)g(Y)]
=
\mathbb E[Xe^{-\lambda X}]
\mathbb E[1-e^{-\lambda X}]
=
\partial_\lambda\psi(u,\lambda)\psi(u,\lambda),
\]
where differentiation under the expectation is justified by
\[
re^{-\lambda r}\le\frac1{e\lambda}.
\]
Combining this with \eqref{eq:pre-burgers},
\[
H(u,\lambda)
\ge
\frac1{10}
\psi(u,\lambda)\partial_\lambda\psi(u,\lambda).
\]
Integrating in \eqref{eq:compensation-integrated} proves
\eqref{eq:burgers-integrated}.
\end{proof}

The next lemma records the characteristic consequence of the Burgers
comparison.

\begin{lemma}[Comparison of the characteristics]
\label{lem:burgers-characteristic}
Fix \(t_0<t_\infty\) and \(\lambda_0>0\), and put
\[
q:=\psi(t_0,\lambda_0).
\]
Define
\[
\lambda(s)
:=
\lambda_0-\frac{q}{10}(s-t_0)
\]
as long as \(\lambda(s)>0\). Then
\[
\psi(s,\lambda(s))\ge q
\]
for every \(s\ge t_0\) such that
\[
s<t_\infty,
\qquad
\lambda(s)>0.
\]
\end{lemma}

\begin{proof}
Fix \(t_1>t_0\) with \(t_1<t_\infty\) and \(\lambda(t_1)>0\), and put
\(v(s):=\psi(s,\lambda(s))\) for \(s\in[t_0,t_1]\).
Lemma~\ref{lem:avalanche-compensation} permits the ordinary chain rule.
Proposition~\ref{prop:burgers-comparison} therefore gives
\[
v'(s)
=\partial_t\psi(s,\lambda(s))
-\frac{q}{10}\,\partial_\lambda\psi(s,\lambda(s))
\ge \frac1{10}\partial_\lambda\psi(s,\lambda(s))(v(s)-q).
\]
Set \(a(s):=\frac1{10}\partial_\lambda\psi(s,\lambda(s))\), which is continuous
on this interval. Multiplication by the integrating factor gives
\[
\frac{d}{ds}\left[
e^{-\int_{t_0}^sa(u)\,du}(v(s)-q)
\right]\ge0.
\]
Since \(v(t_0)=q\), it follows that \(v(s)\ge q\) on \([t_0,t_1]\).
The choice of \(t_1\) was arbitrary.
\end{proof}

\begin{proof}[Proof of Theorem~\ref{thm:sharpness-main}]
Fix \(t<t_\infty\) and \(\lambda>0\), and set
\[
q:=\psi(t,\lambda)>0.
\]
Consider
\[
\lambda(s)
=
\lambda-\frac q{10}(s-t).
\]
Its hitting time of zero is
\[
s_*
=
t+\frac{10\lambda}{q}.
\]
We claim that
\[
s_*\ge t_\infty.
\]

Indeed, if \(s_*<t_\infty\), then
Lemma~\ref{lem:burgers-characteristic} gives
\[
\psi(s,\lambda(s))
\ge q,
\qquad
t\le s<s_*.
\]
Since \(r(s_*)<\infty\) almost surely and \(r(s)\le r(s_*)\),
\[
0
\le
\psi(s,\lambda(s))
\le
\mathbb E[
1-e^{-\lambda(s)r(s_*)}
]
\longrightarrow0
\qquad
(s\uparrow s_*),
\]
by dominated convergence, a contradiction. Therefore
\[
t+\frac{10\lambda}{\psi(t,\lambda)}
\ge t_\infty,
\]
or equivalently
\[
\frac{\psi(t,\lambda)}{\lambda}
\le
\frac{10}{t_\infty-t}.
\]

Since \(r(t)<\infty\) almost surely,
\[
\frac{1-e^{-\lambda r(t)}}{\lambda}
\uparrow r(t)
\qquad
(\lambda\downarrow0).
\]
Monotone convergence yields
\[
R(t)
\le
\frac{10}{t_\infty-t}.
\]
Thus \(R(t)<\infty\) for every \(t<t_\infty\), and hence
\[
t_{\rm sus}\ge t_\infty.
\]
The reverse inequality is immediate. Therefore
\[
t_{\rm sus}=t_\infty=:t_c.
\]
It remains to identify the duration threshold. Since
\(r(t)\le(2d+1)\eta(t)\), finite mean duration implies finite mean range,
so \(t_{\rm d}\le t_{\rm sus}\). Conversely,
Lemma~\ref{lem:range-duration} gives \(D(t)<\infty\) for every
\(t<t_{\rm sus}\), hence \(t_{\rm d}\ge t_{\rm sus}\). Thus
\[
t_{\rm d}=t_{\rm sus}=t_\infty=t_c,
\]
and \eqref{eq:intro-sus-bound} follows.
\end{proof}

\begin{remark}
The quantity
\[
\psi(t,\lambda)
=
\mathbb E[1-e^{-\lambda r(t)}]
\]
may be viewed as the probability that the avalanche range is hit by an
independent Poisson field of intensity \(\lambda\). Infinite
susceptibility corresponds to an infinite right slope at
\(\lambda=0\), whereas an infinite avalanche corresponds to positive
mass at the boundary:
\[
\psi(t,0+)=\theta(t).
\]
The Burgers comparison forces positive values of \(\psi\) to propagate
towards this boundary. Estimate \eqref{eq:intro-sus-bound} expresses
quantitatively that the boundary slope cannot diverge strictly before
positive boundary mass appears.
\end{remark}

\section{Finite-volume amplification}
\label{sec:amplification}

We prove Theorem~\ref{thm:amplification-main}. Fix
\[
d\ge1,
\qquad
t_c<s<t,
\qquad
0<\rho<\frac12,
\]
and put
\[
\Lambda_N:=(\Z/N\Z)^d,
\qquad
V:=N^d,
\qquad
L:=\lfloor\rho N\rfloor.
\]
For \(j\le L\), identify
\[
S_j
:=
\{x\in\Z^d:\|x\|_1=j\}
\]
with its canonical image in \(\Lambda_N\). Since \(2L<N\) for all
sufficiently large \(N\), this identification is unambiguous.

Let
\[
E_N(u)
:=
\{\xi_N(u)\cap S_L\neq\varnothing\}.
\]
The finite- and infinite-volume graphical constructions agree until the
range first reaches \(S_L\). Therefore
\begin{equation}
\label{eq:amp-seed}
\mathbb P(E_N(u))
=
\mathbb P(\xi(u)\cap S_L\neq\varnothing)
\ge
\theta(u).
\end{equation}
In particular,
\[
\mathbb P(E_N(u))
\ge
\theta(s)>0,
\qquad u\ge s.
\]

For \(A\subset\Lambda_N\), write \(a:=|A|\) and
define its translation boundary:
\[
b_A(z)
:=
|\{x\in A:x+z\notin A\}|,
\qquad z\in\Lambda_N.
\]
Since \(A\) and \(A-z\) have equal cardinality,
\begin{equation}
\label{eq:amp-boundary-basic}
b_A(z)
=
\frac12|A\triangle(A-z)|,
\qquad
b_A(-z)=b_A(z),
\qquad
b_A(z+w)\le b_A(z)+b_A(w).
\end{equation}
Moreover,
\begin{equation}
\label{eq:amp-boundary-total}
\sum_{z\in\Lambda_N}b_A(z)
=
a(V-a),
\end{equation}
since both sides count pairs in \(A\times A^c\) according to their
displacement.

\begin{lemma}[Short translations control all translations]
\label{lem:amp-cyclic}
Let
\[
f:\Z/N\Z\to[0,\infty)
\]
satisfy \(f(0)=0\) and
\[
f(r+r')\le f(r)+f(r').
\]
For \(L=\lfloor\rho N\rfloor\), with \(0<\rho<1/2\) fixed, there exists
\(C_\rho<\infty\) such that, for all sufficiently large \(N\),
\[
\sum_{r\in\Z/N\Z}f(r)
\le
C_\rho
\sum_{j=1}^Lf(j).
\]
\end{lemma}

\begin{proof}
Set
\[
S:=\sum_{j=1}^Lf(j).
\]
For \(1\le q<L\), subadditivity gives
\[
f(L)\le f(q)+f(L-q).
\]
Summing over \(q\) yields
\[
(L-1)f(L)\le2S.
\]
Since \(L-1\ge \rho N/2\) for all sufficiently large \(N\), it follows that
\(f(L)\le 4S/(\rho N)\).

For \(0\le r<N\), write \(r=mL+q\), where \(0\le q<L\). Then
\[
f(r)\le mf(L)+f(q).
\]
Also \(L\ge \rho N/2\), hence \(m\le N/L\le2/\rho\), and each value of
\(q\) occurs for at most \(2/\rho+1\) values of \(r\). Therefore
\[
\sum_{r\in\Z/N\Z}f(r)
\le
\frac{2N}{\rho}f(L)
+
\left(\frac2\rho+1\right)S
\le
\left(\frac8{\rho^2}+\frac2\rho+1\right)S.
\]
\end{proof}

Let
\[
H_d
=
(\Z/2\Z)^d\rtimes\mathfrak S_d
\]
be the hyperoctahedral group, acting on \(\Lambda_N\) by sign changes and
coordinate permutations.

\begin{proposition}[Translation expansion]
\label{prop:amp-expansion}
Let \(E\) be an event with positive probability. Suppose that on \(E\),
for every \(1\le j\le L\), we are given a nearest-neighbour edge
\[
(W_j,Z_j),
\qquad
W_j\in S_{j-1},
\quad
Z_j\in S_j,
\]
whose joint conditional law given \(E\) is \(H_d\)-invariant. Define
\[
\kappa
:=
\mathbb E\left[
\mathbf1_E
\sum_{j=1}^L
(\delta_{W_j}+\delta_{Z_j})
\right].
\]
Then
\[
\kappa(\Lambda_N)=2L\mathbb P(E),
\]
and there exists \(c_{d,\rho}>0\) such that, for every
\(A\subset\Lambda_N\), writing \(a=|A|\), and all sufficiently large \(N\),
\begin{equation}
\label{eq:amp-expansion}
\sum_{x\in A,z\in\Lambda_N}
\kappa(z)\mathbf1_{\{x+z\notin A\}}
\ge
c_{d,\rho}\,\kappa(\Lambda_N)
\frac{a(V-a)}V.
\end{equation}
The constant is uniform in \(N\), \(A\), and in the law of the selected
edges.
\end{proposition}

\begin{proof}
Symmetrize the translation boundary by putting
\[
f(z):=\frac1{|H_d|}\sum_{g\in H_d}b_A(gz),
\qquad
I:=\sum_z\kappa(z)b_A(z).
\]
The function \(f\) is subadditive and \(H_d\)-invariant, with
\(f(0)=0\) and \(f(-z)=f(z)\). Invariance of \(\kappa\) and
\eqref{eq:amp-boundary-total} give
\[
I=\sum_z\kappa(z)f(z),
\qquad
\sum_{z\in\Lambda_N}f(z)=a(V-a).
\]

Fix \(z\in S_j\), and let \(a_1,\ldots,a_d\) be its absolute
coordinates, so that \(\sum_i a_i=j\). Permute the coordinates to put
\(a_i\) first, then reflect that coordinate. The difference of the two
resulting vectors is \(2a_i e_1\), and subadditivity and invariance give
\[
f(2a_i e_1)\le2f(z).
\]
Consequently,
\[
f(2j e_1)\le\sum_{i=1}^df(2a_i e_1)\le2d\,f(z),
\]
and the endpoints \(Z_j\) yield
\[
I\ge\frac{\mathbb P(E)}{2d}\sum_{j=1}^Lf(2j e_1).
\]
For each selected edge, \(Z_j-W_j=\pm e_\ell\), so
\[
f(e_1)=f(Z_j-W_j)\le f(Z_j)+f(W_j).
\]
Summing gives \(I\ge\mathbb P(E)Lf(e_1)\).
Together the two estimates imply
\begin{equation}
\label{eq:amp-energy}
I\ge\frac{\mathbb P(E)}{4d}
\left\{
\sum_{j=1}^Lf(2j e_1)+Lf(e_1)
\right\}.
\end{equation}

To control all translations, decompose each vector into its coordinates.
By subadditivity and coordinate symmetry,
\[
\sum_{z\in\Lambda_N}f(z)
\le dN^{d-1}\sum_{r\in\Z/N\Z}f(r e_1).
\]
If \(N\) is odd, multiplication by two permutes \(\Z/N\Z\), so the
last sum equals \(\sum_r f(2r e_1)\). If \(N\) is even, writing each odd
translation as an even translation followed by \(e_1\) gives
\[
\sum_{r\in\Z/N\Z}f(r e_1)
\le
\sum_{r\in\Z/N\Z}f(2r e_1)+\frac N2f(e_1).
\]
Lemma~\ref{lem:amp-cyclic}, applied to \(r\mapsto f(2r e_1)\), gives
\[
\sum_{r\in\Z/N\Z}f(2r e_1)
\le C_\rho\sum_{j=1}^Lf(2j e_1).
\]
Since \(L\ge\rho N/2\) for all sufficiently large \(N\), the extra term
in the even case satisfies
\[
\frac N2f(e_1)\le\frac1\rho Lf(e_1).
\]
Thus in both parity cases we conclude that
\begin{equation}
\label{eq:amp-all-translations}
a(V-a)
\le C_{d,\rho}N^{d-1}
\left\{
\sum_{j=1}^Lf(2j e_1)+Lf(e_1)
\right\}.
\end{equation}
Combining \eqref{eq:amp-energy} and \eqref{eq:amp-all-translations}
yields
\[
I\ge c_{d,\rho}\mathbb P(E)N\frac{a(V-a)}V.
\]
Since
\[
\kappa(\Lambda_N)=2L\mathbb P(E)\le2\rho N\mathbb P(E),
\]
the preceding estimate, after changing the constant, gives
\[
I\ge c_{d,\rho}\,\kappa(\Lambda_N)\frac{a(V-a)}V.
\]
Finally,
\[
I=\sum_{x\in A,z\in\Lambda_N}
\kappa(z)\mathbf1_{\{x+z\notin A\}},
\]
which proves \eqref{eq:amp-expansion}.
\end{proof}

\begin{remark}
We need both endpoints of the selected edges to handle parity.
The vertices \(Z_j\) alone control the even translations \(2j e_1\).
When \(N\) is even, these do not control odd translations. Using also
\(W_j\), the nearest-neighbour relation \(Z_j-W_j=\pm e_\ell\) yields
\(f(e_1)\le f(Z_j)+f(W_j)\), which bridges the two parity classes.
For odd \(N\), multiplication by \(2\) is invertible modulo \(N\), and
this extra ingredient is unnecessary.
\end{remark}

\begin{remark}
\label{rem:poincare}
Proposition~\ref{prop:amp-expansion} also gives a spectral-gap bound for
the auxiliary translation walk. Put \(m:=\kappa(\Lambda_N)\) and
\(\mu:=\kappa/m\), and let \(P_\mu(x,y):=\mu(y-x)\). By symmetry of
\(\kappa\), this kernel is reversible for the uniform law. The corresponding
continuous-time generator is
\[
L_\kappa f(x)
:=
\sum_z\kappa(z)\bigl(f(x+z)-f(x)\bigr)
=
m(P_\mu-I)f(x).
\]
The proposition gives the conductance bound
\[
\inf_{0<|A|\le V/2}
\frac{\sum_z\mu(z)b_A(z)}{|A|}
\ge\frac{c_{d,\rho}}2.
\]
Hence Cheeger's inequality \cite[Theorem~13.10]{LevinPeresWilmer2017}
implies that the spectral gap \(\lambda_\kappa\) of \(-L_\kappa\) satisfies
\[
\lambda_\kappa
=
m\,\operatorname{gap}(I-P_\mu)
\ge\frac{c_{d,\rho}^2}{8}\,m.
\]
Conversely, a spectral-gap estimate of the form
\(\lambda_\kappa\ge c'_{d,\rho}\kappa(\Lambda_N)\)
automatically provides an alternative proof of
Proposition~\ref{prop:amp-expansion}. Indeed, writing \(\pi\) for the
uniform law, the Dirichlet form associated with \(L_\kappa\) is
\[
\mathcal E_\kappa(f,f)
:=
-\langle f,L_\kappa f\rangle_\pi
=
\frac1{2V}\sum_{x,z}\kappa(z)
\bigl(f(x+z)-f(x)\bigr)^2.
\]
Thus
\[
\mathcal E_\kappa(\mathbf1_A,\mathbf1_A)
=
\frac1V\sum_z\kappa(z)b_A(z),
\qquad
\operatorname{Var}_\pi(\mathbf1_A)
=
\frac{|A|(V-|A|)}{V^2}.
\]
Poincar\'e's inequality
\cite[Lemma~13.7 and Remark~13.8]{LevinPeresWilmer2017},
\[
\lambda_\kappa\operatorname{Var}_\pi(f)
\le \mathcal E_\kappa(f,f),
\]
then implies \eqref{eq:amp-expansion}, up to constants.
\end{remark}

We now apply Proposition~\ref{prop:amp-expansion} to the avalanche.
On \(E_N(u)\), connectedness of \(\xi_N(u)\) implies that, for every
\(1\le j\le L\), the range contains an edge
\[
(W_j,Z_j),
\qquad
W_j\in S_{j-1},
\quad
Z_j\in S_j.
\]
Choose one uniformly among all such edges, independently for each \(j\)
conditional on the range. This selection is equivariant under \(H_d\), so
its conditional law given \(E_N(u)\) is \(H_d\)-invariant. Set
\[
\kappa_u
:=
\mathbb E\left[
\mathbf1_{E_N(u)}
\sum_{j=1}^L
(\delta_{W_j}+\delta_{Z_j})
\right].
\]
Since
\[
\kappa_u(\Lambda_N)=2L\mathbb P(E_N(u)),
\]
while \eqref{eq:amp-seed} gives
\[
\mathbb P(E_N(u))\ge\theta(s),
\qquad
u\ge s,
\]
and, since \(L\ge \rho N/2\) for all sufficiently large \(N\),
\[
\kappa_u(\Lambda_N)\ge \rho\,\theta(s)N.
\]
Hence Proposition~\ref{prop:amp-expansion} gives
\(c_{d,s,\rho}>0\) such that, uniformly in \(u\ge s\),
\(A\subset\Lambda_N\), and large \(N\),
\begin{equation}
\label{eq:amp-avalanche-expansion}
\sum_{x\in A,z}
\kappa_u(z)\mathbf1_{\{x+z\notin A\}}
\ge
c_{d,s,\rho}N
\frac{|A|(V-|A|)}V.
\end{equation}

By Proposition~\ref{prop:hierarchical-decomposition}, conditionally on
the current range \(A\) at level \(u\), every \(x\in A\) carries marks at
rate one, and a mark rooted at \(x\) grafts an independent
\(u\)-avalanche translated by \(x\).

For such a graft, the event \(E_N(u)\) and the selected endpoints refer
to the independent avalanche. Let \(D_A(x)\) be the number of distinct
sites outside \(A\) among the selected translated endpoints, with
\(D_A(x)=0\) when \(E_N(u)\) fails. A vertex on \(S_j\) can occur only as
\(Z_j\), and possibly also as \(W_{j+1}\); hence every selected endpoint
has multiplicity at most two. Therefore
\[
0\le D_A(x)\le2L,
\qquad
D_A(x)
\ge
\frac12\mathbf1_{E_N(u)}
\sum_{j=1}^L
\left(
\mathbf1_{\{x+W_j\notin A\}}
+
\mathbf1_{\{x+Z_j\notin A\}}
\right).
\]
Summing over the possible roots, taking expectations, and using the
definition of \(\kappa_u\) and \eqref{eq:amp-avalanche-expansion} gives,
after decreasing the constant,
\begin{equation}
\label{eq:amp-graft-mean}
\sum_{x\in A}\mathbb E[D_A(x)]
\ge c_{d,s,\rho}N\frac{|A|(V-|A|)}V.
\end{equation}
The expectation is over the independent graft and the auxiliary selection
of its crossing edges, with \(A\) fixed. The actual number of new sites
created by a graft rooted at \(x\) dominates \(D_A(x)\).

\begin{proposition}[Supercritical covering]
\label{prop:amp-spanning}
For every \(t_c<s<t\) and every fixed \(0<\rho<1/2\),
\[
\mathbb P(
\xi_N(t)=\Lambda_N
\mid E_N(s)
)
\longrightarrow1.
\]
Consequently,
\[
\liminf_{N\to\infty}P_N(t)
\ge
\theta(s)>0.
\]
\end{proposition}

\begin{proof}
Condition throughout on \(E_N(s)\). Connectedness gives
\(r_N(s)\ge L+1\), and the range can only increase. For a nonempty
\(A\subset\Lambda_N\), put
\[
\Phi(A):=\frac V{|A|}-1.
\]
Let \(a:=|A|\ge L\). A graft rooted at \(x\in A\) adds
\(\Delta:=|(x+\xi_N'(u))\setminus A|\) new sites, where
\(\xi_N'(u)\) is an independent \(u\)-avalanche. Since
\(\Delta\ge D_A(x)\), \(D_A(x)\le2L\le2a\), and
\(v\mapsto v/(a+v)\) is increasing,
\[
\Phi(A)-\Phi(A\cup(x+\xi_N'(u)))
=\frac{V\Delta}{a(a+\Delta)}
\ge\frac{VD_A(x)}{a(a+D_A(x))}
\ge\frac{VD_A(x)}{3a^2}.
\]
Let \(\mathcal L_u\) denote the generator of the range process.
Each root has rate one, so \eqref{eq:amp-graft-mean} gives, uniformly
for \(u\ge s\) and \(a\ge L\),
\[
-\mathcal L_u\Phi(A)
\ge\frac V{3a^2}\sum_{x\in A}\mathbb E[D_A(x)]
\ge\frac{c_{d,s,\rho}N}{3}\frac{V-a}{a}
=\frac{c_{d,s,\rho}N}{3}\Phi(A).
\]
The conditioning event belongs to \(\mathcal H_s\), and all subsequent
grafts retain their independent laws. Dynkin's formula and Gronwall's
inequality therefore yield
\begin{equation}
\label{eq:amp-lyapunov}
\mathbb E[\Phi(\xi_N(t))\mid E_N(s)]
\le\frac VL
\exp\left\{-\frac{c_{d,s,\rho}N}{3}(t-s)\right\}.
\end{equation}
If \(A\neq\Lambda_N\), then \(\Phi(A)\ge1/(V-1)\). Consequently,
\begin{equation}
\label{eq:amp-cover-error}
\mathbb P(\xi_N(t)\neq\Lambda_N\mid E_N(s))
\le\frac{V(V-1)}L
\exp\left\{-\frac{c_{d,s,\rho}N}{3}(t-s)\right\}
\longrightarrow0.
\end{equation}
Finally,
\[
P_N(t)\ge
\mathbb P(E_N(s))
\mathbb P(\xi_N(t)=\Lambda_N\mid E_N(s)),
\]
and \eqref{eq:amp-seed} gives
\(\liminf_{N\to\infty}P_N(t)\ge\theta(s)>0\).
\end{proof}

\begin{proof}[Proof of Theorem~\ref{thm:amplification-main}]
This is Proposition~\ref{prop:amp-spanning}.
\end{proof}

\begin{remark}
The estimate \eqref{eq:amp-cover-error} controls covering throughout the
growth from a macroscopic arm to the whole torus. Since \(V=N^d\) and \(L\ge\rho N/2\) for all sufficiently large \(N\),
its prefactor satisfies
\[
\frac{V(V-1)}L\le\frac2\rho N^{2d-1}.
\]
Thus the comparison also supplies a covering bound after a sprinkling
\(C(\log N)/N\), with \(C\) sufficiently large depending on
\(d,s,\rho\). 
\end{remark}

\section{Subcritical locality}
\label{sec:subcritical-moments}

We prove Theorem~\ref{thm:subcritical-main}. Recall that
\[
R(t)=\mathbb E r(t),
\qquad
t_c=t_{\rm sus},
\]
so that \(R(t)<\infty\) for every \(t<t_c\). The hierarchical decomposition
first upgrades this finite susceptibility to exponential integrability of
the avalanche range. Local injectivity then compares infinite- and
finite-volume ranges, while Lemma~\ref{lem:range-duration} transfers the
resulting convergence of mean ranges to mean durations. The following proof is inspired by \cite{riordan2015evolution}.

\begin{proposition}[Subcritical exponential moments]
\label{prop:subcritical-moments}
For every \(t<t_c\), there exists \(\lambda_t>0\) such that
\[
\mathbb E[e^{\lambda_t r(t)}]<\infty.
\]
In particular, \(\mathbb E[r(t)^p]<\infty\) for every integer \(p\ge1\).
\end{proposition}

\begin{proof}
Fix \(L\ge1\) and, for \(\lambda\ge0\), set
\[
F_L(u,\lambda):=\mathbb E[e^{\lambda(r(u)\wedge L)}].
\]
By Proposition~\ref{prop:hierarchical-decomposition}, conditionally on a
current range \(A\) of size \(x\), marks occur at total rate \(x\), and
each mark grafts an independent avalanche \(B\) with
\(|B|\stackrel{\mathrm d}=r(u)\). Since
\[
|A\cup(z+B)|\le x+|B|,
\]
if \(x<L\), then
\[
e^{\lambda((x+|B|)\wedge L)}-e^{\lambda x}
\le
e^{\lambda x}\bigl(e^{\lambda(|B|\wedge L)}-1\bigr),
\]
whereas the truncated observable is constant when \(x\ge L\).
Conditioning on the current range and using the independence of the graft
therefore gives, for almost every \(u<t_c\),
\begin{equation}
\label{eq:exp-transport}
\partial_uF_L(u,\lambda)
\le
\bigl(F_L(u,\lambda)-1\bigr)
\partial_\lambda F_L(u,\lambda).
\end{equation}
Indeed, the first factor produced by the rate is bounded by
\[
\mathbb E\!\left[
r(u)e^{\lambda r(u)}\mathbf1_{\{r(u)<L\}}
\right]
\le \partial_\lambda F_L(u,\lambda).
\]

We claim that exponential integrability propagates over every interval
\([s,s+h]\) such that \(hR(s)<1\). Suppose
\(\mathbb E[e^{\alpha_0r(s)}]<\infty\) for some \(\alpha_0>0\).
Since
\[
\frac{\mathbb E[e^{\alpha r(s)}]-1}{\alpha}\longrightarrow R(s)
\qquad(\alpha\downarrow0),
\]
we may choose \(0<\alpha<\alpha_0\) such that, writing
\(q:=\mathbb E[e^{\alpha r(s)}]\),
\[
h(q-1)<\alpha.
\]
For \(s\le u\le s+h\), set
\[
\lambda(u):=\alpha-(q-1)(u-s),
\qquad
G_L(u):=F_L(u,\lambda(u)).
\]
Then \(\lambda(u)>0\), and \eqref{eq:exp-transport} yields
\[
G_L'(u)
\le
\partial_\lambda F_L(u,\lambda(u))
\bigl(G_L(u)-q\bigr),
\qquad
G_L(s)\le q.
\]
Hence \(G_L(u)\le q\) on \([s,s+h]\). In particular,
\[
F_L(s+h,\alpha-h(q-1))\le q,
\qquad
\alpha-h(q-1)>0.
\]
Letting \(L\to\infty\) gives
\[
\mathbb E[e^{(\alpha-h(q-1))r(s+h)}]\le q<\infty.
\]

Now fix \(t<t_c\). Since \(R(t)<\infty\), choose a partition
\(0=s_0<\cdots<s_m=t\) with mesh smaller than \(R(t)^{-1}\).
The initial range \(r(0)=|\mathcal N_d|=2d+1\) is deterministic.
If \(r(s_j)\) has an exponential moment, then, by monotonicity of \(R\),
\[
(s_{j+1}-s_j)R(s_j)
\le (s_{j+1}-s_j)R(t)<1,
\]
so the preceding propagation argument gives an exponential moment for
\(r(s_{j+1})\). Iterating finitely many times yields
\[
\mathbb E[e^{\lambda_t r(t)}]<\infty
\]
for some \(\lambda_t>0\).
\end{proof}

The preceding argument is independent of the geometry of \(\Z^d\).
Geometry enters only when these range moments are used to compare the infinite
lattice with the torus.

\begin{lemma}[Local injectivity]
\label{lem:local-injectivity}
Let
\[
\varpi_N:\Z^d\to\Lambda_N
\]
be the canonical projection. If \(A\subset\Z^d\) is connected and
\[
|A|<N,
\]
then \(\varpi_N\) is injective on \(A\).

Consequently, for every \(K\) with
\[
K+2d<N,
\]
the infinite- and finite-volume avalanche constructions can be coupled
in physical update time, using the same refresh variables, so as to agree,
including their durations, until their common range first contains at
least \(K\) sites.
\end{lemma}

\begin{proof}
A connected set with \(k\) vertices has graph diameter at most \(k-1\).
Thus, if \(|A|<N\),
\[
\operatorname{diam}_{\ell^1}(A)
\le
|A|-1
<
N.
\]
On the other hand, distinct \(x,y\in\Z^d\) satisfying
\[
\varpi_N(x)=\varpi_N(y)
\]
have
\[
\|x-y\|_1\ge N.
\]
Hence \(\varpi_N\) is injective on \(A\).

For the second assertion, couple the elementary Bak--Sneppen updates
in \(\Z^d\) and on the torus using the same refresh variables, as long as
the explored range projects injectively. Avalanche ranges are connected, and one elementary
Bak--Sneppen update can add at most \(2d\) new sites to the range. Hence,
while the current range contains fewer than \(K\) sites, after the next
update it contains fewer than \(K+2d<N\) sites and still projects
injectively. The two explorations therefore remain identical until the
claimed stopping time.
\end{proof}

\begin{proposition}[Subcritical finite-volume approximation]
\label{prop:subcritical-finite-volume}
For every \(t<t_c\),
\begin{equation}
\label{eq:range-finite-volume}
\sup_{0\le u\le t}
|R_N(u)-R(u)|
\longrightarrow0.
\end{equation}
Moreover,
\begin{equation}
\label{eq:duration-finite-volume}
D_N(t)\longrightarrow D(t).
\end{equation}
Finally, for every fixed \(K\), for all sufficiently large \(N\) and every
\(1\le k\le K\),
\begin{equation}
\label{eq:local-duration-identification}
\mathbb E\left[
\eta_N(t)\mathbf1_{\{r_N(t)=k\}}
\right]
=
\mathbb E\left[
\eta(t)\mathbf1_{\{r(t)=k\}}
\right].
\end{equation}
\end{proposition}

\begin{proof}
Take
\[
K_N:=\lfloor N/3\rfloor.
\]
For large \(N\),
\[
K_N+2d<N,
\]
so Lemma~\ref{lem:local-injectivity} couples the two avalanches identically
until their common range reaches \(K_N\). In particular, on
\(\{r(u)<K_N\}\),
\[
r_N(u)=r(u).
\]
Since
\[
r_N(u)\le N^d,
\]
we obtain
\[
|R_N(u)-R(u)|
\le
N^d
\mathbb P(r(u)\ge K_N)
+
\mathbb E[
r(u)
\mathbf1_{\{r(u)\ge K_N\}}
].
\]

Fix an integer \(p>d\). By monotonicity in the level parameter and
Proposition~\ref{prop:subcritical-moments},
\[
\mathbb P(r(u)\ge K_N)
\le
\frac{\mathbb E[r(t)^p]}{K_N^p},
\]
and
\[
\mathbb E[
r(u)
\mathbf1_{\{r(u)\ge K_N\}}
]
\le
\frac{\mathbb E[r(t)^p]}{K_N^{p-1}},
\]
uniformly for \(u\le t\). Hence
\[
\sup_{u\le t}|R_N(u)-R(u)|
\le
C_{p,t}
\left(
N^{d-p}+N^{1-p}
\right)
\longrightarrow0,
\]
which proves \eqref{eq:range-finite-volume}.

By Lemma~\ref{lem:range-duration},
\[
D_N(t)
=
\exp\left\{\int_0^tR_N(u)\,du\right\},
\qquad
D(t)
=
\exp\left\{\int_0^tR(u)\,du\right\}.
\]
Therefore \eqref{eq:range-finite-volume} implies
\eqref{eq:duration-finite-volume}.

Finally fix \(K\), and take \(N\) large enough that
\[
K+1+2d<N.
\]
Apply Lemma~\ref{lem:local-injectivity} with \(K+1\). The two avalanches
then agree until their common range first contains at least \(K+1\) sites.
On the event that the final range has size \(k\le K\), this stopping level
is never reached. Hence the entire finite-volume avalanche is identical to
the infinite-volume avalanche, including its duration, and conversely. This
proves \eqref{eq:local-duration-identification}.
\end{proof}

\begin{proof}[Proof of Theorem~\ref{thm:subcritical-main}]
The exponential-moment statement is Proposition~\ref{prop:subcritical-moments},
and the convergence of the mean ranges is \eqref{eq:range-finite-volume}.
\end{proof}

\section{From avalanches to the stationary state}
\label{sec:stationary}

We now complete the passage from avalanche criticality to the stationary
state. The argument is the finite-volume mechanism introduced by Meester
and Znamenski \cite{MeesterZnamenski2003,MeesterZnamenski2004}, but we give
it in full because it is dimension-free and is an essential part of the
present theorem. Write
\[
V:=|\Lambda_N|=N^d,
\qquad
\nu_h:=\mathcal L(h+E),\quad E\sim\operatorname{Exp}(1).
\]
All neighbourhoods on the torus are understood with repetitions removed.

\begin{proposition}[Finite-volume ergodicity]
\label{prop:finite-volume-ergodicity}
For every \(N\), the Bak--Sneppen chain on \(\Lambda_N\) has a unique
stationary law \(\pi_N\). From every deterministic initial configuration,
its law converges to \(\pi_N\) in total variation. Moreover \(\pi_N\) is
translation invariant.
\end{proposition}

\begin{proof}
If one update refreshes the whole torus, then after one step the
configuration consists of independent \(\operatorname{Exp}(1)\) variables,
independently of the initial state, and the conclusion is immediate.
Assume henceforth that one update does not cover \(\Lambda_N\).

Fix \(h>0\), write \(P\) for the transition kernel and \(V=N^d\), and set
\[
C_h:=\{f:f_x>h\text{ for every }x\in\Lambda_N\}.
\]
With probability at least \(q_h:=e^{-(2d+1)hV}>0\), every fitness sampled
in the next \(V\) updates exceeds \(h\), uniformly in the current state.
On this event no value at most \(h\) is created. As long as such values
remain, the selected minimum is at most \(h\) and is refreshed, so their
number decreases at each update. After at most \(V\) updates none remains.
Hence
\begin{equation}
\label{eq:uniform-hit-Ch}
P^V(f,C_h)\ge q_h,
\qquad f\in[0,\infty)^{\Lambda_N}.
\end{equation}
Iteration gives a geometric tail for the hitting time of \(C_h\),
uniformly in the initial state.

For the minorization on \(C_h\), choose a finite nearest-neighbour walk
starting at the first updated site and visiting every vertex. At each
intermediate update, require that the prescribed next vertex receive a
fitness at most \(h\), while every other refreshed vertex receives a
fitness above \(h\); at the final update require all fresh values to
exceed \(h\). At each intermediate step the prescribed next vertex is
then the unique site with fitness at most \(h\), so the avalanche follows
the walk, covers the torus, and terminates. These finitely many
requirements have positive probability, so \(P_N(h)>0\). Since the
avalanche duration is integer-valued, there exists \(m\ge1\) such that
\[
p_h:=\mathbb P\bigl(r_N(h)=V,\ \eta_N(h)=m\bigr)>0.
\]
This range--duration event is measurable with respect to the censored
level-\(h\) exploration, whose law, up to translation of its starting
point, does not depend on \(f\in C_h\). On this event every site has
been refreshed, so Proposition~\ref{prop:hierarchical-decomposition}
implies that the terminal fitnesses are independent with common law
\(\nu_h\), independently of the censored exploration. Thus, setting
\[
\mu_{N,h}:=\nu_h^{\otimes\Lambda_N},
\]
we obtain
\begin{equation}
\label{eq:small-set-minorization}
P^m(f,\cdot)\ge p_h\,\mu_{N,h}(\cdot),
\qquad f\in C_h.
\end{equation}

Combining \eqref{eq:uniform-hit-Ch} and
\eqref{eq:small-set-minorization} by the Markov property gives
\begin{equation}
\label{eq:global-doeblin}
P^{V+m}(f,\cdot)\ge q_hp_h\,\mu_{N,h}(\cdot)
\qquad\text{for every }f.
\end{equation}
Hence \(\lambda\mapsto\lambda P^{V+m}\) is a contraction in total
variation with factor \(1-q_hp_h\). The space of probability measures is
complete for this metric, so Banach's fixed-point theorem gives a unique
fixed point \(\pi_N\). Since \(\pi_NP\) is also
\(P^{V+m}\)-invariant, uniqueness gives \(\pi_NP=\pi_N\), so \(\pi_N\)
is the unique invariant law of \(P\). Iterating the contraction and using
that Markov kernels are nonexpansive in total variation yields
\[
\|\lambda P^n-\pi_N\|_{\rm TV}
\le
(1-q_hp_h)^{\lfloor n/(V+m)\rfloor}
\|\lambda-\pi_N\|_{\rm TV}
\longrightarrow0
\]
for every initial law \(\lambda\).

Finally, start the chain from independent \(\operatorname{Exp}(1)\)
fitnesses. Its law is translation invariant, ties have probability zero at
every time, and the Bak--Sneppen dynamics preserves translation invariance.
The total-variation convergence just proved therefore implies that its limit
\(\pi_N\) is translation invariant.
\end{proof}

We next introduce the locking variables. Start the finite-volume chain from
independent \(\operatorname{Exp}(1)\) fitnesses and set
\(Y_{N,x}(0)=0\). If \(i_N(n)\) is the site of minimum fitness at time
\(n\) and
\[
M_N(n):=F_{N,i_N(n)}(n),
\]
define
\begin{equation}
\label{eq:locking-recursion}
Y_{N,x}(n+1)
=
\begin{cases}
0,&x\in i_N(n)+\mathcal N_d,\\
\max\{Y_{N,x}(n),M_N(n)\},&x\notin i_N(n)+\mathcal N_d.
\end{cases}
\end{equation}
For \(y,z\ge0\), set
\[
H_y(z):=
\begin{cases}
0,&z<y,\\
1-e^{-(z-y)},&z\ge y.
\end{cases}
\]

\begin{proposition}[Conditional exponential representation]
\label{prop:conditional-exponential}
Conditionally on \(Y_N(n)=(Y_{N,x}(n))_{x\in\Lambda_N}\), the fitnesses
\((F_{N,x}(n))_{x\in\Lambda_N}\) are independent and
\[
F_{N,x}(n)\sim \nu_{Y_{N,x}(n)}.
\]
Consequently, for distinct \(x_1,\ldots,x_k\in\Lambda_N\) and
\(z_1,\ldots,z_k\ge0\),
\begin{equation}
\label{eq:locking-exp-representation}
\mathbb P\bigl(F_{N,x_j}(n)\le z_j,\ 1\le j\le k\bigr)
=
\mathbb E\left[
\prod_{j=1}^k H_{Y_{N,x_j}(n)}(z_j)
\right].
\end{equation}
\end{proposition}

\begin{proof}
The assertion holds at time zero. Assume it holds at time \(n\), and
condition on \(Y_N(n)\), the location \(i_N(n)\), and the value \(M_N(n)\)
of the minimum. By the lack-of-memory property, the unrefreshed coordinates
remain independent and can be written as
\[
\max\{Y_{N,x}(n),M_N(n)\}+E'_{N,x},
\]
where the \(E'_{N,x}\) are independent \(\operatorname{Exp}(1)\) variables.
The refreshed coordinates are fresh independent exponentials. Thus,
conditionally on these variables, the fitnesses at time \(n+1\) are
independent exponentials above the thresholds in
\eqref{eq:locking-recursion}. Since this conditional law depends on the
conditioning variables only through \(Y_N(n+1)\), the same statement holds
conditionally on \(Y_N(n+1)\) alone. This proves the claim by induction.
\end{proof}

Put
\[
G_N(n,t):=\mathbb P(Y_{N,0}(n)\le t).
\]
Translation invariance of the initial law and of the dynamics makes the
choice of the origin immaterial.

\begin{proposition}[Subcritical locking thresholds]
\label{prop:locking-subcritical}
For every \(t<t_c\),
\begin{equation}
\label{eq:locking-below}
\lim_{N\to\infty}\limsup_{n\to\infty}G_N(n,t)=0.
\end{equation}
\end{proposition}

\begin{proof}
Let
\[
K_N(n,t):=\sum_{x\in\Lambda_N}
\mathbf1_{\{Y_{N,x}(n)\le t\}}.
\]
Let \(C_t\) be the set of configurations with all fitnesses above \(t\).
By \eqref{eq:uniform-hit-Ch}, its hitting time has finite mean. Discard the
initial segment up to the first visit to \(C_t\), and from then on decompose
the trajectory into successive \(t\)-avalanches. If \(T_j\) is the starting
time of the \(j\)-th one and \((r_j,\eta_j)\) its range--duration pair, then
\(T_{j+1}-T_j=\eta_j\), and these pairs are i.i.d. with the law of
\((r_N(t),\eta_N(t))\). At the first update of the avalanche the current
minimum exceeds \(t\), so every site outside its eventual range has its
locking threshold raised above \(t\). During the remaining updates the
selected minima are at most \(t\), and only sites in the avalanche range
are refreshed. Hence
\[
K_N(n,t)\le r_j,
\qquad T_j<n\le T_{j+1}.
\]
For all sufficiently large \(N\), a one-step noncovering \(t\)-avalanche
has positive probability, so the cycle length has span one. The delayed
renewal--reward theorem and translation invariance therefore give
\[
\limsup_{n\to\infty}G_N(n,t)
\le
\frac{\mathbb E[\eta_N(t)r_N(t)]}{V D_N(t)}.
\]
For every fixed \(K\ge1\), using \(r_N(t)\le V\),
\begin{equation}
\label{eq:locking-subcritical-estimate}
\limsup_{n\to\infty}G_N(n,t)
\le
\frac K V
+
\frac{
\mathbb E[\eta_N(t)\mathbf1_{\{r_N(t)>K\}}]
}{D_N(t)}.
\end{equation}
Now
\[
\mathbb E[\eta_N(t)\mathbf1_{\{r_N(t)>K\}}]
=
D_N(t)-
\sum_{k=1}^K
\mathbb E[\eta_N(t)\mathbf1_{\{r_N(t)=k\}}].
\]
By Proposition~\ref{prop:subcritical-finite-volume}, as \(N\to\infty\)
this converges to
\(\mathbb E[\eta(t)\mathbf1_{\{r(t)>K\}}]\), while
\(D_N(t)\to D(t)<\infty\). Letting first \(N\to\infty\) and then
\(K\to\infty\) in \eqref{eq:locking-subcritical-estimate} proves
\eqref{eq:locking-below}.
\end{proof}

\begin{proposition}[Supercritical locking thresholds]
\label{prop:locking-supercritical}
For every \(b>t_c\),
\begin{equation}
\label{eq:locking-above}
\lim_{N\to\infty}\liminf_{n\to\infty}G_N(n,b)=1.
\end{equation}
\end{proposition}

\begin{proof}
Fix \(t_c<t<b\). Starting from the process defining \(G_N\), ignore the
initial segment up to the first visit to \(C_t\). Thereafter run successive
\(t\)-avalanches until the first one covering \(\Lambda_N\), and then
continue until the next visit to \(C_b\).

The initial delay and the final recovery have finite mean by
\eqref{eq:uniform-hit-Ch}. Moreover \(P_N(t)>0\) by
Proposition~\ref{prop:finite-volume-ergodicity}, and
\(D_N(t)<\infty\) by Lemma~\ref{lem:range-duration}, so the waiting time
for the first covering \(t\)-avalanche also has finite mean.

At the end of this covering avalanche, every site has been refreshed.
Conditionally on its \(t\)-censored exploration, the terminal fitnesses
are independent with law \(\nu_t\), by
Proposition~\ref{prop:hierarchical-decomposition}. Censor the subsequent
recovery at level \(b\): reveal fitnesses at most \(b\) exactly and reveal
larger fitnesses only through their exceedance of \(b\). At the designated
return to \(C_b\), each terminal value has been observed only to exceed
\(b\). Thus, conditionally on the censored cycle history, the terminal
vector has law
\[
\mu_{N,b}=\nu_b^{\otimes\Lambda_N}.
\]
Treat the trajectory up to the first such return as the initial delay,
and repeat the construction from that point. Fresh future refresh
variables and this censored regeneration make the successive pairs of
waiting and recovery lengths independent and identically distributed.
Let \(A_N\) be the time spent in successive \(t\)-avalanches up to and
including the first covering one in a cycle started with law
\(\mu_{N,b}\), and let \(B_N\) be the subsequent recovery time to \(C_b\).
For large \(N\), the cycle length has span one, since a
one-step noncovering \(t\)-avalanche may be inserted before the first
covering one with positive probability.

At the end of the covering avalanche all locking thresholds are at most
\(t\), and throughout the recovery every selected minimum is at most \(b\).
Thus
\[
K_N(n,b)=V
\]
during the whole recovery phase. By the delayed renewal--reward theorem,
\begin{equation}
\label{eq:supercritical-reward}
\liminf_{n\to\infty}G_N(n,b)
\ge
\frac{\mathbb E B_N}{\mathbb E A_N+\mathbb E B_N}.
\end{equation}

Since the successive \(t\)-avalanches are i.i.d. until the first covering
one, Wald's identity gives
\begin{equation}
\label{eq:covering-wald}
\mathbb E A_N=\frac{D_N(t)}{P_N(t)}.
\end{equation}

To bound \(B_N\) from below, use the decorated construction of
Proposition~\ref{prop:hierarchical-decomposition}. Since the covering
\(t\)-avalanche has range \(\Lambda_N\), every site is active in the level
extension from \(t\) to \(b\). For each \(x\in\Lambda_N\), let \(\tau_x\)
be the first mark above \(t\) on the vertical line through \(x\). Then
\(\tau_x-t\) are i.i.d. \(\operatorname{Exp}(1)\), and, conditionally on
\(\tau_x=u\le b\), the corresponding graft contributes a duration of mean
\(D_N(u)\). Hence, by Tonelli and monotonicity of \(D_N\),
\begin{equation}
\label{eq:recovery-lower-bound}
\mathbb E B_N
\ge
V\int_t^b e^{-(u-t)}D_N(u)\,du
\ge
V(1-e^{-(b-t)})D_N(t).
\end{equation}

Combining
\eqref{eq:supercritical-reward}--\eqref{eq:recovery-lower-bound},
\begin{equation}
\label{eq:locking-supercritical-estimate}
1-\liminf_{n\to\infty}G_N(n,b)
\le
\frac{1}{V(1-e^{-(b-t)})P_N(t)}.
\end{equation}
Choose \(s\in(t_c,t)\). By Theorem~\ref{thm:amplification-main},
\[
\liminf_{N\to\infty}P_N(t)\ge\theta(s)>0.
\]
Since \(V=N^d\to\infty\), the right-hand side of
\eqref{eq:locking-supercritical-estimate} tends to zero, proving
\eqref{eq:locking-above}.
\end{proof}

\begin{proposition}[Locking-threshold concentration and product limit]
\label{prop:locking-product-limit}
For every \(\varepsilon>0\),
\begin{equation}
\label{eq:locking-concentration}
\lim_{N\to\infty}\limsup_{n\to\infty}
\mathbb P\bigl(|Y_{N,0}(n)-t_c|>\varepsilon\bigr)=0.
\end{equation}
Consequently, for every fixed collection of distinct sites
\(x_1,\ldots,x_k\in\Z^d\), the corresponding marginal of
\(\pi_N\) converges to \(\nu_{t_c}^{\otimes k}\).
\end{proposition}

\begin{proof}
For \(0<\varepsilon<t_c\), Propositions~\ref{prop:locking-subcritical}
and \ref{prop:locking-supercritical}, applied at \(t_c-\varepsilon\) and
\(t_c+\varepsilon\), give \eqref{eq:locking-concentration}; for
\(\varepsilon\ge t_c\), the lower-tail event is empty and the same
conclusion follows from the supercritical bound. For fixed
\(x_1,\ldots,x_k\), and \(N\) large enough that their images in
\(\Lambda_N\) are distinct, translation invariance and the union bound imply
\[
\lim_{N\to\infty}\limsup_{n\to\infty}
\mathbb P\left(
\max_{1\le j\le k}|Y_{N,x_j}(n)-t_c|>\varepsilon
\right)=0.
\]

Fix \(z_1,\ldots,z_k\ge0\) and define
\[
\Psi(y_1,\ldots,y_k):=
\prod_{j=1}^k H_{y_j}(z_j).
\]
The function \(\Psi\) is bounded and continuous. If
\[
\omega_\varepsilon:=
\sup_{\max_j|y_j-t_c|\le\varepsilon}
|\Psi(y_1,\ldots,y_k)-\Psi(t_c,\ldots,t_c)|,
\]
then \(\omega_\varepsilon\to0\) as \(\varepsilon\downarrow0\). By
Proposition~\ref{prop:conditional-exponential},
\[
\begin{aligned}
&\left|
\mathbb P(F_{N,x_j}(n)\le z_j,\ 1\le j\le k)
-\prod_{j=1}^k H_{t_c}(z_j)
\right|\\
&\qquad\le
\omega_\varepsilon
+
\mathbb P\left(
\max_{1\le j\le k}|Y_{N,x_j}(n)-t_c|>\varepsilon
\right).
\end{aligned}
\]
For fixed \(N\), Proposition~\ref{prop:finite-volume-ergodicity} lets
\(n\to\infty\) and replaces the law of the chain by \(\pi_N\). Letting
next \(N\to\infty\) and finally \(\varepsilon\downarrow0\) proves
convergence of the joint distribution functions to
\(\prod_j H_{t_c}(z_j)\), the distribution function of
\(\nu_{t_c}^{\otimes k}\).
\end{proof}

\begin{proof}[Proof of Theorem~\ref{thm:stationary-main}]
Proposition~\ref{prop:locking-product-limit} gives the asserted convergence
in exponential coordinates. Under the inverse transformation
\(u=1-e^{-f}\), a random variable with law
\(t_c+\operatorname{Exp}(1)\) becomes uniform on
\([1-e^{-t_c},1]\). This proves the theorem.
\end{proof}

\appendix
\section{Use of AI}
\label{appendix-AI}

The (compensated) Laplace transform of the nonlinear additive process in which \(X_t\) is replaced by \(X_t+X_t'\), with \(X_t'\) an independent copy of \(X_t\), at rate \(X_t\), satisfies an inviscid Burgers equation. Starting from \(X_0=1\), its solution is given by the Borel(\(t\)) distribution, the law of the total progeny of a Poisson(\(t\)) Galton--Watson tree. The phase transition at \(t=1\) can then be identified by the method of characteristics; see \cite{MenonPego2004,DeaconuTanre2000}.

This mechanism turns out to be closely related to the branching-process description of Bak--Sneppen avalanches in the mean-field setting \cite{deBoerDerrida1994}, and for two years, we had been looking for a way to relate it to the evolution of the avalanche range in our spatial setting. During that period, we sought a suitable comparison or domination by solutions of the Burgers equation, both through our own efforts and by prompting earlier AI models. These attempts were unsuccessful and, in particular, repeatedly led to incorrect boundary-layer arguments by AI.

In August 2026, ChatGPT 5.6 Sol produced the key argument used in this paper: restrict attention to grafts involving large avalanches, derive a first-moment lower bound on the number of newly created vertices, and combine this estimate with concavity to obtain the Burgers-type comparison (Lemma~\ref{lem:translation-creation} and Proposition~\ref{prop:burgers-comparison}).

The amplification argument was subsequently obtained by the same model through a spectral-gap approach, along the lines of Remark~\ref{rem:poincare}. We then suggested looking for a purely combinatorial proof, leading to the argument retained in the paper.

Finally, finiteness of all polynomial moments in the subcritical regime was first established by the model through a recursive argument. We then suggested seeking a recursive proof of exponential moments inspired by \cite{riordan2015evolution}; the resulting argument was developed with AI assistance.

Claude did the simulation leading to the  picture.

\bibliographystyle{amsplain}
\bibliography{references}

\end{document}